\documentclass{amsart}
\usepackage{graphicx}
\usepackage{amscd}
\usepackage{amsmath}
\usepackage{amsfonts}
\usepackage{amssymb}
\theoremstyle{plain}

\newtheorem{corollary}{Corollary}

\newtheorem{definition}{Definition}
\newtheorem{example}{Example}

\newtheorem{lemma}{Lemma}

\newtheorem{remark}{Remark}

\newtheorem{theorem}{Theorem}
\numberwithin{equation}{section}

\begin{document}
\title[Classification of Banach Spaces]{Topological and Cardinality Properties of Certain Sets of Classes of Banach
Spaces }
\author{Eugene Tokarev}
\address{B.E. Ukrecolan, 33-81 Iskrinskaya str., 61005, Kharkiv-5, Ukraine}
\email{tokarev@univer.kharkov.ua}
\thanks{This paper is in final form and no version of it will be submitted for
publication elsewhere.}
\dedicatory{Dedicated to the memory of S. Banach.}
\begin{abstract}Classes of Banach spaces that are finitely, strongly finitely or elementary
equivalent are introduced. On sets of these classes topologies are defined in
such a way that sets of defined classes become compact totally disconnected
topological spaces. Results are used in the problem of synthesis of Banach
spaces, and to describe omittable spaces that are defined below.
\end{abstract}
\tableofcontents
\maketitle

A general classification scheme, suggested by author for Banach spaces
includes their partition by dimension (which may be regarded as horizontal
strips of the proposed coordinate system) and a partition by the relation
$\sim_{f}$ of finite equivalence (defined below), that consists of classes of
finitely equivalent Banach spaces (and which may be regarded as vertical
strips of the mentioned generalized coordinate system). That class of finite
equivalence, which contains a given Banach space $X$ will be denoted by
$X^{f}$. These classes may be further participated in a few ways. In the paper
we consider two of such partition: by the relation $\approx_{\phi}$ of strong
finite equivalence (that consists of classes $X^{\phi}$ of strong finite
equivalence) and by the relation of elementary equivalence $\equiv$, that
generates a partition of the class $\mathcal{B}$ of all Banach spaces in
classes $X^{\xi}$ of pairwice elementary equivalent Banach spaces.

Certainly, the partition by dimension enjoys such nice property as
well-ordering. If we fix a cardinal $\varkappa$, then the set $\{\mathcal{B}%
_{\tau}:\tau\leq\varkappa\}$ may be equipped with the order topology, under
which it became compact.

Although other mentioned partitions are not well-ordered (in the partition
into classes of elementary equivalence no order will be defined at all), on
these 'coordinate strips' suitable topologies that turn them in topological
compact spaces may be defined as well.

Dependence of properties of Banach spaces on their position in the classes
$X^{f}$, $X^{\phi}$ and $X^{\xi}$ is studied in other papers of author. Here
only cardinal and topological properties of the classification will be studied.

\section{Classification of Banach spaces by dimension}

\subsection{Definitions}

\textit{Ordinals} will be denoted by small Greece letters $\alpha$, $\beta$,
$\gamma$. \textit{Cardinals} are identified with the least ordinals of given
cardinality and are denoted either by $\iota$, $\tau$, $\varkappa$, $\varrho$,
or by using Hebrew letter $\gimel$ (may be with indices). As usual, $\omega$
and $\omega_{1}$ denote respectively the first infinite and the first
uncountable cardinal (= ordinal).

For a cardinal $\tau$ its \textit{predecessor} (i.e. the least cardinal,
strongly greater then $\tau$) is denoted by $\tau^{+}$.

\textit{The} \textit{confinality of} $\tau,$ $\operatorname{cf}(\tau)$ is the
least cardinality of a set $A\subset\tau$ such that $\tau=\sup A$.

Let $A$, $B$ be sets. The symbol $^{B}A$ denotes the set of all functions from
$B$ to $A$.

In a general case, the cardinality of the set $^{B}A$ is denoted either by
$\operatorname{card}(A)^{\operatorname{card}(B)}$ or by $\varkappa^{\tau}$, if
$\operatorname{card}(A)=\varkappa$; $\operatorname{card}(B)=\tau$.

For an ordinal $\gamma$ the cardinal $\gimel(\gamma)$ is given by induction:%
\[
\gimel(0)=\omega;\text{ \ }\gimel(\alpha+1)=2^{\gimel(\alpha)};\text{
\ }\gimel(\gamma)=\cup\{\gimel(\alpha):\alpha<\gamma\},
\]
if $\gamma$ is a limit ordinal.

The symbol $\exp(\tau)$ (or, equivalently, $2^{\tau}$) denotes the cardinality
of the set $\operatorname{Pow}\left(  \tau\right)  $ of all subsets of $\tau$.

The cardinal $\exp(\omega)$ (i.e., the cardinality of continuum) will be also
denoted by $\frak{c}$.

Let $X$ be a Banach space. Below all Banach spaces will be considered over the
field $\mathbb{R}$ of real scalars.

Any set of elements $\{w_{\beta}:\beta<\alpha\}$\ (of arbitrary nature),
indexed by ordinals $<\alpha$ will be called an $\alpha$-sequence.

An $\alpha$-sequence $\{x_{\beta}:\beta<\alpha\}$ of elements of $X$ is said
to be

\begin{itemize}
\item \textit{Spreading,} if for any $n<\omega$, any $\varepsilon>0$, any
scalars $\{a_{k}:k<n\}$ and any choosing of $i_{0}<i_{1}<...<i_{n-1}<\alpha$;
$j_{0}<j_{1}<...<j_{n-1}<\alpha$ the following equality holds:
\[
\left\|  \sum\nolimits_{k<n}a_{k}x_{i_{k}}\right\|  =\left\|  \sum
\nolimits_{k<n}a_{k}x_{j_{k}}\right\|  ;
\]

\item \textit{Symmetric,} if for any $n<\omega$, any finite subset
$I\subset\mathbb{\alpha}$ of cardinality $n$, any rearrangement $\varsigma$ of
elements of $I$ and any scalars $\{a_{i}:i\in I\}$,
\[
\left\|  \sum\nolimits_{i\in I}a_{i}z_{i}\right\|  =\left\|  \sum
\nolimits_{i\in I}a_{\varsigma(i)}z_{i}\right\|  .
\]

\item \textit{Subsymmetric}, if it is both spreading and 1-unconditional.
\end{itemize}

Let $C<\infty$ be a constant. Two $\alpha$-sequences $\{x_{\beta}:\beta
<\alpha\}$ and $\{y_{\beta}:\beta<\alpha\}$ are said to be $C$%
\textit{-equivalent} if for any finite subset $I=\{i_{0}<i_{1}<...<i_{n-1}\}$
of $\alpha$ and for any choosing of scalars $\{a_{k}:k<n\}$
\[
C^{-1}\left\|  \sum\nolimits_{k<n}a_{k}x_{i_{k}}\right\|  \leq\left\|
\sum\nolimits_{k<n}a_{k}y_{i_{k}}\right\|  \leq C\left\|  \sum\nolimits_{k<n}%
a_{k}x_{i_{k}}\right\|  .
\]

Two $\alpha$-sequences $\{x_{\beta}:\beta<\alpha\}$ and $\{y_{\beta}%
:\beta<\alpha\}$ are said to be equivalent if they are $C$-equivalent for some
$C<\infty$.

\subsection{Banach spaces spanned by subsymmetric sequences}

It may be shown that if $\frak{x}$ is not symmetric (resp., is not equivalent
to a symmetric sequence) then the cardinality $\operatorname{card}%
(\left\lfloor \frak{x}\right\rfloor \cap\mathcal{B}_{\varkappa})=2^{\varkappa
}$ (resp., the cardinality $\operatorname{card}\left(  (\left\lfloor
\frak{x}\right\rfloor \cap\mathcal{B}_{\varkappa})^{\approx}\right)
=2^{\varkappa}$).

Let $\varkappa$ be a cardinal; $\sigma:\varkappa\rightarrow\varkappa$ be a
transposition (i.e. one-to-one mapping of $\varkappa$ onto $\varkappa$).

It will be said that $\sigma$ is almost identical if the correlation
\[
\gamma_{1}<\gamma_{2}\Rightarrow\sigma\gamma_{1}<\sigma\gamma_{2}\text{, where
}\gamma_{1}\text{, }\gamma_{2}<\varkappa
\]
may get broken at most finitely many times.

\begin{lemma}
Let $\frak{x}=\{x_{n}:n<\omega\}$ be a spreading sequence that is not symmetric.

Let $\sigma:\omega\rightarrow\omega$ be a transposition, which is not almost identical.

If sequences $\{x_{n}:n<\omega\}$ and $\{x_{\sigma n}:n<\omega\}$ are
equivalent then both of them are equivalent to a symmetric sequences.
\end{lemma}

\begin{proof}
Consider a sequence $\{y_{\alpha}:\alpha<\omega^{2}\}$, which is given by
\[
\{y_{\omega\cdot2k+n}=x_{n};\text{ \ }y_{\omega\cdot\left(  2k+1\right)
+n}=x_{\sigma n}:k<\omega;\text{ }n<\omega\}.
\]
This sequence is equivalent to a sequence $\{y_{\alpha}^{\prime}:\alpha
<\omega^{2}\}$ that belongs to the tower $\left\lfloor \frak{x}\right\rfloor
$, generated by $\frak{x}=\{x_{n}:n<\omega\}$.

Certainly, this is equivalent to
\[
C^{-1}\left\|  \sum\nolimits_{k<n}a_{k}x_{i_{k}}\right\|  <\left\|
\sum\nolimits_{k<n}a_{k}y_{\alpha_{k}}\right\|  <C\left\|  \sum\nolimits_{k<n}%
a_{k}x_{i_{k}}\right\|
\]

for every $n<\omega$; every scalars $\left(  a_{k}\right)  _{k<n}$ every
choice $\alpha_{0}<\alpha_{1}<...<\alpha_{n-1}<\omega^{2}$, $i_{0}%
<i_{1}<...<i_{n-1}<\omega$ and some $C<\infty$ that depends only on
equivalence constant between $\{x_{n}:n<\omega\}$ and $\{x_{\sigma n}%
:n<\omega\}$. Since $\sigma$ has only a finite number of inversions, our
definition of $\{y_{\alpha}:\alpha<\omega^{2}\}$ implies that $\frak{x}$ is
equivalent to a symmetric sequence
\end{proof}

\begin{theorem}
Let $\alpha$, $\beta$ be ordinals; $\omega\leq\beta\leq\alpha$. Let
$\{x_{\gamma}:\gamma<\alpha\}$ be a subsymmetric $\alpha$-sequence which is
not equivalent to any $\alpha$-symmetric sequence. Let $X_{\alpha
}=\operatorname{span}(\{x_{\gamma}:\gamma<\alpha\})$; $X_{\beta}%
=\operatorname{span}(\{x_{\gamma}:\gamma<\beta\})$. If $\beta^{\omega}<\alpha$
then spaces $X_{\alpha}$ and $X_{\beta}$ are not isomorphic.
\end{theorem}

\begin{proof}
Assume that $X_{\alpha}$ is isomorphic to a subspace $Z$\ of $X_{\beta}$.

Let $I:X_{\alpha}\rightarrow X_{\beta}$ be the corresponding operator of
isomorphic embedding. Without loss of generality it may be assumed that an
image of element $x_{\gamma}$ ($\gamma<\alpha$) in $X_{\beta}$ is a finite
linear combination with rational coefficients of some $x_{\zeta}$'s
($\zeta<\beta$):
\[
Ix_{\gamma}=\sum\nolimits_{k=0}^{n(\gamma)}a_{k}^{\gamma}x_{\zeta_{k}\left(
\gamma\right)  };\text{ \ }\zeta_{0}\left(  \gamma\right)  <\zeta_{1}\left(
\gamma\right)  <...<\zeta_{n(\gamma)}\left(  \gamma\right)  <\beta.
\]
Thus, to any $x_{\gamma}$ corresponds a finite sequence of rational numbers
$\left(  a_{k}^{\gamma}\right)  _{k=0}^{n(\gamma)}$.

Let $\left(  p_{n}\right)  _{n<\omega}$ be a numeration of all finite
sequences of rationals.

The set $\left(  p_{n}\right)  _{n<\omega}$ generates a partition of $\alpha$
on parts $\left(  P_{n}\right)  _{n<\omega}$ in a following way: an ordinal
$\gamma<\alpha$ belongs to $P_{n}$ if $\left(  a_{k}^{\gamma}\right)
_{k=0}^{n(\gamma)}=p_{n}$.

Since $\alpha>\beta^{\omega}\geq\omega^{\omega}$, one of $P_{n}$'s should
contain a sequence $\{\gamma_{i}:i<\delta\}$ of ordinals, which order type
$\delta$ (in a natural order : $i<j$ implies that $\gamma_{i}<\gamma_{j}$) is
greater then $\beta^{\omega}$. Clearly, for all such $\gamma_{i}$,
\[
Ix_{\gamma_{i}}=\sum\nolimits_{k=0}^{n}a_{k}x_{\zeta_{k}\left(  \gamma
_{i}\right)  },
\]
where $n$ and $\left(  a_{k}\right)  _{k=1}^{n}$ do not depend on $i$.

A set of all sequences $\zeta_{0}\left(  \gamma_{i}\right)  <\zeta_{1}\left(
\gamma_{i}\right)  <...<\zeta_{n}\left(  \gamma_{i}\right)  <\beta$ cannot be
ordered to have the order type $>\beta^{n}$. Hence, conditions $\gamma_{i_{1}%
}<\gamma_{i_{2}}\Rightarrow\zeta_{0}\left(  \gamma_{i_{1}}\right)  <\zeta
_{0}\left(  \gamma_{i_{2}}\right)  $ must get broken for infinite many pairs
$\gamma_{i_{1}}$, $\gamma_{i_{2}}$. The inequality
\[
C^{-1}\left\|  \sum\nolimits_{k<m}b_{k}x_{\gamma_{k}}\right\|  \leq\left\|
\sum\nolimits_{k<m}b_{k}\left(  \sum\nolimits_{k=0}^{n}a_{k}x_{\zeta
_{k}\left(  \gamma_{i}\right)  }\right)  \right\|  \leq C\left\|
\sum\nolimits_{k<m}b_{k}x_{\gamma_{k}}\right\|  ,
\]
where $C=d(X_{\alpha},Z)$, shows that $\{x_{\gamma}:\gamma<\alpha\}$ is
equivalent to a symmetric sequence.
\end{proof}

\begin{theorem}
Let $\frak{x}=\{x_{n}:n<\omega\}$ be a subsymmetric sequence, which is not
isomorphic to a symmetric one. Then for any cardinal $\varkappa\geq\omega$
there exists $2^{\varkappa}$ pairwice non isomorphic Banach spaces of
dimension $\varkappa$, which belongs to the same tower $\left\lfloor
\frak{x}\right\rfloor $.
\end{theorem}

\begin{proof}
Let $\{x_{\gamma}:\gamma<\varkappa\}$ be a subsymmetric $\varkappa$-sequence.
Let $I=\left\langle I,\ll\right\rangle $ be a linearly ordered set of
cardinality $\varkappa$; $J=\left\langle I,<^{\prime}\right\rangle $ --
another linear ordering on $I$.

Consider families $\{x_{i}:i\in I\}$ and $\{x_{j}:j\in J\}$ that are indexed
(and ordered) by elements of $I$ and $J$ respectively.

Let $X_{I}=\operatorname{span}\{x_{i}:i\in I\}$; $X_{J}=\operatorname{span}%
\{x_{j}:j\in J\}$.

Certainly, $X_{I}$ and $X_{J}$ are isomorphic if only if there are one-to-one
mappings of embedding $u:I\rightarrow J$ and $w:J\rightarrow I$, which are
almost monotone in a following sense:
\begin{align*}
i_{1} &  \ll i_{2}\Rightarrow u\left(  i_{1}\right)  <^{\prime}u\left(
i_{i}\right)  \text{ \ for all but finitely many pairs }i_{1},i_{2}\in I;\\
j_{1} &  <^{\prime}j_{2}\Rightarrow w\left(  j_{1}\right)  \ll w\left(
j_{i}\right)  \text{ \ for all but finitely many pairs }j_{1},j_{2}\in J.
\end{align*}
Since there exists $2^{\varkappa}$ orderings of $I$ for any pair of which such
almost monotone mappings do not exist, this prove the theorem.
\end{proof}

\subsection{Classes $\mathcal{B}_{\varkappa}$}

Let $\mathcal{B}$ be a proper class of all Banach spaces (isometric spaces are
identified). For any finite or infinite cardinal number $\varkappa$ let
$\mathcal{B}_{\varkappa}$ be a set of all Banach spaces of dimension
$\varkappa$:%
\[
\mathcal{B}_{\varkappa}=\{X\in\mathcal{B}:\dim X=\varkappa\}
\]

Recall that a \textit{dimension} $\dim(X)$ of a Banach space $X$ is the least
cardinality of a subset $A\subset X$, which \textit{linear span}
$\operatorname{lin}(A)$ is dense in $X$ (equivalently: such $A\subset X$ that
a \textit{closure} $\overline{\operatorname{lin}(A)}$, which will be in the
future denoted by $\operatorname{span}(A)$, is the whole space $X$). If $X$ is
of infinite dimension then its dimension $\dim(X)$ is exactly equal to its
\textit{density character} $\operatorname{dens}(X)$ -- the least cardinality
of a subset $B\subset X$, which is dense in $X$.

This partition may be regarded as an analogue of horizontal strips of some
generalized coordinate system for Banach spaces.

Consider a question on cardinality of classes $\mathcal{B}_{\varkappa}$.

Immediately $\operatorname*{card}\mathcal{B}_{0}=\operatorname*{card}%
\mathcal{B}_{1}=1$.

Indeed, the inclusion $X\in\mathcal{B}_{0}$ means that $X$ is $0$-dimensional,
i.e., $X=\{0\}$.

Every $1$-dimensional Banach space $X=\operatorname*{span}e$ may be identified
with the scalar field $\mathbb{R}$ with the norm $\left\|  \lambda e\right\|
=\left|  \lambda\right|  $.

\begin{theorem}
For every finite cardinal $n\geq2$ the cardinality $\operatorname*{card}%
\mathcal{B}_{n}=\exp\left(  \omega\right)  =\frak{c}$.
\end{theorem}

\begin{proof}
Let $n<\omega$. Consider the vector space $\mathbb{R}^{n}$, equipped with a
Hausdorff topology (since all such topologies on finite-dimensional vector
space $\mathbb{R}^{n}$ are equivalent, one may assume, e.g. that
$\mathbb{R}^{n}$ is equipped with the Euclidean metric).

To any $n$-dimensional Banach space $X$ corresponds its unit ball $B(X)=\{x\in
X:\left\|  x\right\|  \leq1\}$, which is the central symmetric closed body in
$\mathbb{R}^{n}$. Because of the set of such bodies in $\mathbb{R}^{n}$ is at
most of cardinality $\frak{c}$, it is clear that $\operatorname*{card}%
\mathcal{B}_{n}\leq\exp\left(  \omega\right)  $.

To show the converse inequality, consider a collections of spaces
$l_{p}^{\left(  n\right)  }$ ($1\leq p\leq\infty$).\ Recall that the upper
index $^{\left(  n\right)  }$ pointed out on dimension. Surely, all these
spaces are pairwice non-isometric and, hence
\[
\operatorname*{card}\mathcal{B}_{n}=\exp\left(  \omega\right)  =\frak{c.}%
\]
\end{proof}

Notice that all Banach spaces of given finite dimension are pairwice
isomorphic (= linearly homeomorphic). This is not a case when
infinite-dimensional Banach spaces are under consideration.

E.g., the same argument shows that $\operatorname*{card}\mathcal{B}_{\omega
}=\exp\left(  \omega\right)  $. However, because of all spaces $l_{p}$ (
$1\leq p<\infty$) are of dimension $\omega$ and are pairwice non-isomorphic,
the more powerful result is valid. To formulate it we introduce some notations.

Let $X$ be a Banach space. Let $X^{\approx}$ denotes a set of all Banach
spaces $Y$ that are isomorphic to $X$ (shortly: $Y\approx X$ and, thus,
$X^{\approx}=\{Y\in\mathcal{B}:Y\approx X\}$). Surely, all spaces from
$X^{\approx}$ are of the sane dimension as $X$. Obviously, $\approx$ is an
equivalence relation. which participates $\mathcal{B}_{\varkappa}$ into parts
modulo $\approx$. Put%
\[
\mathcal{B}_{\varkappa}^{\approx}=\{X^{\approx}:X\in\mathcal{B}_{\varkappa
}\}.
\]

Obviously, $\operatorname*{card}\mathcal{B}_{n}^{\approx}=1$ for all
$n<\omega$. Since $l_{p}$ and $l_{q}$ are isomorphic if only if $p=q$, it
follows that $\operatorname*{card}\mathcal{B}_{\omega}^{\approx}=\exp\left(
\omega\right)  $. Indeed,
\[
\exp\left(  \omega\right)  =\operatorname*{card}\mathbb{R=}%
\operatorname*{card}\{l_{p}:1\leq p<\infty\}\leq\operatorname*{card}%
\mathcal{B}_{\omega}^{\approx}\leq\operatorname*{card}\mathcal{B}_{\omega
}=\exp\left(  \omega\right)  .
\]

The similar result is valid for every infinite cardinal.

\begin{theorem}
Let $\varkappa>\omega$ be a cardinal; $\{x_{\alpha}:\alpha<\varkappa\}$ and
$\{y_{\alpha}:\alpha<\varkappa\}$ be subsymmetric sequences;
$X=\operatorname{span}\{x_{\alpha}:\alpha<\varkappa\}$; $Y=\operatorname{span}%
\{y_{\alpha}:\alpha<\varkappa\}$. If spaces $X$ and $Y$ are isomorphic then
$\varkappa$-sequences $\{x_{\alpha}:\alpha<\varkappa\}$ and $\{y_{\alpha
}:\alpha<\varkappa\}$ are equivalent.
\end{theorem}

\begin{proof}
Let $u:X\rightarrow Y$ be an isomorphism; $\left\|  u\right\|  \left\|
u^{-1}\right\|  =c$. It may be assumed that $ux_{\alpha}\in Y$ is represented
as a block
\[
ux_{\alpha}=\sum\nolimits_{k=1}^{n(\alpha)}a_{k}^{\alpha}y_{\beta_{k}\left(
\alpha\right)  }%
\]
for some sequence of rational scalars $\left(  a_{k}^{\alpha}\right)  _{k\leq
n(\alpha)}$ and finite $n(\alpha)$. Moreover, it may be assumed that this
blocks are not intersected, i.e., that any member of a given sequence $\left(
\beta_{k}\left(  \alpha\right)  \right)  _{k\leq n(\alpha)}$ belongs only to
this block.

Since $\varkappa$ is uncountable, among such blocks there is an infinite
number of identical ones that differs only in sequences $\left(  \beta
_{k}\left(  \alpha\right)  \right)  _{k\leq n(\alpha)}$. Let $A\subset
\varkappa$ be such that all elements $\{ux_{\alpha}:\alpha\in A\}$ are
represented by identical blocks:
\[
ux_{\alpha}=\sum\nolimits_{k=1}^{n}a_{k}y_{\beta_{k}\left(  \alpha\right)
}\text{ \ for }\alpha\in A.
\]

Let $\left(  b_{i}\right)  $ be a sequence of scalars; $a=\max_{1\leq k\leq
n}\left(  \left|  a_{k}\right|  \right)  $.

Then for any finite subset $A^{\prime}\subset A$, $A^{\prime}=\{\alpha
_{i}\}_{i=1}^{m}$, because of unconditionality of sequences $\left(
x_{\alpha}\right)  $ and $\left(  y_{\alpha}\right)  $,
\begin{align*}
\left\|  \sum\nolimits_{i=1}^{m}b_{i}x_{\alpha_{i}}\right\|   &  \geq
c\left\|  \sum\nolimits_{i=1}^{m}b_{i}\left(  \sum\nolimits_{k=1}^{n}%
a_{k}y_{\beta_{k}\left(  i\right)  }\right)  \right\|  \\
&  \geq c\left\|  \sum\nolimits_{i=1}^{m}b_{i}nay_{\beta_{1}\left(  i\right)
}\right\|  \geq cna\left\|  \sum\nolimits_{i=1}^{m}b_{i}y_{\beta_{1}\left(
i\right)  }\right\|  .
\end{align*}

Analogously,\ holds the converse inequality that proves the theorem.
\end{proof}

\begin{corollary}
The cardinality of the set of all Banach spaces of dimension $\varkappa$ is
\[
\operatorname{card}\left(  \mathcal{B}_{\varkappa}^{\approx}\right)
=\operatorname{card}\left(  \mathcal{B}_{\varkappa}\right)  =\exp\left(
\omega\right)  .
\]
\end{corollary}

\begin{proof}
The inequality $\operatorname{card}\left(  \mathcal{B}_{\varkappa}^{\approx
}\right)  \geq2^{\varkappa}$ follows from the previous theorem and results of
the previous subsection. The inverse inequality is obvious.
\end{proof}

\begin{remark}
It is of interest that different sets $\mathcal{B}_{\varkappa}$ and
$\mathcal{B}_{\tau}$ may be of the same cardinality. The appearance of a such
case depends on the model of the set theory that we chose as the base of all
functional analysis.

E.g., if one assume the Martin axiom $MA$ with the negation of continuum
hypothesis $\urcorner CH$ ($MA+\urcorner CH$) then for all cardinals
$\varkappa$ such that $\omega<\varkappa<2^{\omega}$
\[
\operatorname{card}\mathcal{B}_{\varkappa}=\operatorname{card}\mathcal{B}%
_{\omega}=\exp\left(  \omega\right)  .
\]
\end{remark}

As it was noted in the introduction, for a given cardinal $\tau$ one may
define a set$\frak{B}^{\left(  \tau\right)  }$, which elements are sets
$\mathcal{B}_{\varkappa}$ for all cardinals $\varkappa\leq\tau$, namely:%
\[
\frak{B}^{\left(  \tau\right)  }=\{\mathcal{B}_{\varkappa}:\varkappa\leq
\tau;\text{ \ }\varkappa\text{ is a cardinal}\}.
\]

This set may be identified with the set of all ordinals $\{\alpha:\alpha
\leq\tau\}$ and, as it is well known, may be endowed with the order topology
(say, $\mathcal{T}_{\left(  o\right)  }$), under which it becomes a compact
topological space. This result we formulate as follows.

\begin{theorem}
Let $\frak{B}^{\left(  \tau\right)  }$ be defined as above. Then $\left\langle
\frak{B}^{\left(  \tau\right)  },\mathcal{T}_{\left(  o\right)  }\right\rangle
$ is a compact topological space for every cardinal $\tau$.
\end{theorem}

\section{Classification of Banach spaces by elementary equivalence}

\subsection{Properties of ultrapowers}

The notion of ultraproducts and ultrapowers of Banach spaces was introduced by
D. Dacunha-Castelle and J.-L. Krivine \cite{DCCK} as an analogue of the
model-theoretical notion of ultraproducts and ultrapowers of models for the
first-order language (more details see in the next section).

General properties of ultraproducts and ultrapowers were luxuriously exposed
in the S. Heinrich's paper \cite{Henr}. Here will be presented results that
are not contained in \cite{Henr}.

\begin{definition}
Let $I$ be a set; $\operatorname{Pow}(I)$ be a set of all its subsets. An
ultrafilter $D$ over $I$ is a subset of $\operatorname{Pow}(I)$ with following properties:

\begin{itemize}
\item $I\in D$;

\item  If $A\in D$ and $A\subset B\subset I$, then $B\in I$;

\item  If $A$, $B$ $\in D$ then $A\cap B\in D$;

\item  If $A\in D$, then $I\backslash A\notin D$.
\end{itemize}
\end{definition}

\begin{definition}
Let $I$ be a set; $D$ be an ultrafilter over $I$; $\{X_{i}:i\in I\}$ be a
family of Banach spaces. An \textit{ultraproduct } $(X_{i})_{D}$ is given by a
quotient space
\[
(X)_{D}=l_{\infty}\left(  X_{i},I\right)  /N\left(  X_{i},D\right)  ,
\]
where $l_{\infty}\left(  X_{i},I\right)  $ is a Banach space of all families
$\frak{x}=\{x_{i}\in X_{i}:i\in I\}$, for which
\[
\left\|  \frak{x}\right\|  =\sup\{\left\|  x_{i}\right\|  _{X_{i}}:i\in
I\}<\infty;
\]
$N\left(  X_{i},D\right)  $ is a subspace of $l_{\infty}\left(  X_{i}%
,I\right)  $, which consists of such $\frak{x}$'s$\ $that
\[
\lim_{D}\left\|  x_{i}\right\|  _{X_{i}}=0.
\]
\end{definition}

If all $X_{i}$'s are all equal to a space $X\in\mathcal{B}$ then an
ultraproduct is said to be an \textit{ultrapower} and is denoted by $\left(
X\right)  _{D}$.

An operator $d_{X}:X\rightarrow\left(  X\right)  _{D}$ that asserts to any
$x\in X$ an element $\left(  x\right)  _{D}\in\left(  X\right)  _{D}$, which
is generated by a stationary family $\{x_{i}=x:i\in I\}$ is called the
\textit{canonical embedding }of $X$ into its ultrapower $\left(  X\right)
_{D}$.

\begin{definition}
Let $D$ be an ultrafilter, $\varkappa$ be a cardinal. An ultrafilter $D$ is called

\begin{itemize}
\item $\varkappa$-regular if there exists a subset $G\subset D$ of cardinality
$\operatorname{card}(G)=\varkappa$ such that any $i\in I$ belongs only to a
finite number of sets $e\in G$;

\item $\varkappa$-complete if an intersection of a nonempty subset $G\subset
D$ of cardinality $<\varkappa$ belongs to $D$, i.e., if
\[
G\subset D\text{ \ and \ }\operatorname{card}(G)<\varkappa\text{ \ implies
that }\cap G\in D;
\]

\item  Principal (or non-free), if there exists $i\in I$ such that
$D=\{e\subset I:i\in e\}$ (such ultrafilter is said to be generated by $i\in
I$);

\item  Free, if it is not a principal one;

\item  Countably incomplete if it is not $\omega_{1}$-complete.
\end{itemize}
\end{definition}

\begin{theorem}
Let $X$ be a Banach space, $\dim(X)=\varkappa$; $D$ be an ultrafilter.

The canonical embedding $d_{X}:X\rightarrow\left(  X\right)  _{D}$ maps $X$
\textbf{on }$\left(  X\right)  _{D}$ if and only if $D$ is $\varkappa^{+}$-complete.
\end{theorem}

\begin{proof}
Let $D$ be $\varkappa^{+}$-complete. If $\operatorname{card}(I)\leq\varkappa$
then $D$ is non-free and, hence, $d_{X}X=\left(  X\right)  _{D}$. Assume that
$\operatorname{card}(I)>\varkappa$ and that $\left(  x_{i}\right)  _{D}%
\in\left(  X\right)  _{D}\backslash d_{X}X$. Of course, $\left(  x_{i}\right)
_{D}$ is generated by a such family $\left(  x_{i}\right)  _{i\in I}$ that for
any $\varepsilon>0$
\[
I_{0}(\varepsilon)=\{i\in I:\left\|  x_{i}-x\right\|  >\varepsilon\}\in D
\]
for any $x\in X$. Consider a function $F:I_{0}\rightarrow X$, which is given
by $f(i)=x_{i}$. Since $\dim(X)=\varkappa$, a partition
\[
I=\cup\{f^{-1}\left(  x\right)  :x\in I_{0}\left(  \varepsilon\right)
\}\cup\{I\backslash I_{0}\left(  \varepsilon\right)  \}
\]
participate $I$ in $\tau$ sets where $\tau<\varkappa^{+}$. Since $D$ is
$\varkappa^{+}$-complete, one of sets of the partition belongs to $D$. Because
of $I_{0}\left(  \varepsilon\right)  \in D$, $I\backslash I_{0}\left(
\varepsilon\right)  \notin D$. Hence there exists such $x\in X$ that
$f^{-1}\left(  x\right)  \in D$. Clearly, this contradicts with the assumption
$\left(  x_{i}\right)  _{D}\in\left(  X\right)  _{D}\backslash d_{X}X$.

Conversely, assume that $d_{X}X=\left(  X\right)  _{D}$. Since $\dim
(X)=\varkappa$ there exists $\varepsilon>0$ and a set $A\subset X$ such that
$\operatorname{card}(A)=\varkappa$ and $\left\|  a-b\right\|  \geq\varepsilon$
for all $a\neq b\in A$.

Let $I=\cup\{I_{a}:a\in B\}$ be a partition of $I$ in $\tau
=\operatorname{card}(B)<\varkappa^{+}$ parts. Let a function $F:I\rightarrow
A$ be given by:
\[
F(i)=a\text{ \ if and only if }i\in I_{a}\text{.}%
\]
Then $\left(  F\left(  i\right)  \right)  _{D}\in\left(  X\right)  _{D}%
=d_{X}X$, and, hence, $\left(  F\left(  i\right)  \right)  _{D}=d_{X}\left(
a\right)  $ for some $a\in A$. Because of $\left\|  a-b\right\|
\geq\varepsilon$ for all $a\neq b\in A$, $F^{-1}(a)\in D$. However, by our
assumption, $F^{-1}(a)=I_{a}$. Thus, $I_{a}\in D$ and, since a partition of
$I$ by $\tau$ parts was arbitrary, $D$ is $\varkappa^{+}$-complete.
\end{proof}

\begin{corollary}
Let $X$ be a Banach space, $D$ be an ultrafilter.

\begin{enumerate}
\item  If $\dim(X)<\omega$ then $\left(  X\right)  _{D}=X$;

\item  If $\omega\leq\dim(X)<\upsilon$ where $\upsilon$ is the first
measurable cardinal, and $D$ is a free ultrafilter then $\left(  X\right)
_{D}\backslash d_{X}X\neq\varnothing$;

\item  If $X$ is of infinite dimension and ultrafilter $D$ is countably
incomplete, then $\left(  X\right)  _{D}\backslash d_{X}X\neq\varnothing$ too.
\end{enumerate}
\end{corollary}

\begin{proof}
Since any ultrafilter is $\omega$-complete (by definitions), all results
follows from the preceding theorem.
\end{proof}

Nevertheless, by choosing an ultrafilter, a difference $\left(  X\right)
_{D}\backslash d_{X}X$ may be maiden arbitrary large.

\begin{theorem}
Let $X\in\mathcal{B}_{\varkappa}$; $D$ be a $\tau$-regular ultrafilter over a
set $I$ of cardinality $\operatorname{card}(I)=\tau$. Then
\[
\dim((X)_{D})=(\dim(X))^{\tau}=\varkappa^{\tau}.
\]
\end{theorem}

\begin{proof}
Let $A^{\prime}\subset X$ be a set of cardinality $\varkappa$ which is dense
in $X$. Immediately,
\[
\dim\left(  \left(  X\right)  _{D}\right)  \leq\dim\left(  l_{\infty}\left(
X,I\right)  \right)  \leq\operatorname{card}(^{I}A)=\left(
\operatorname{card}\left(  A\right)  \right)  ^{\operatorname{card}\left(
I\right)  }=\varkappa^{\tau}.
\]

Assume now that $A\subset A^{\prime}$ is a set of the same cardinality
$\varkappa$ such that for some $\varepsilon>0$ and any $a,b\in A$, $a\neq b$,
$\left\|  a-b\right\|  \geq\varepsilon$. Let $\{a_{i}:i<\varkappa\}$ be a
numeration of elements of $A$. Any finite subset $s=\{x_{0},...,x_{n-1}%
\}\subset A$ spans a finite dimensional subspace $X_{s}$ of $X$.

A set $B$ of all subspaces of kind $X_{s}$ of $X$ is also of cardinality
$\varkappa$.

Let $\{X_{\alpha}:\alpha<\varkappa\}$ be a numeration of $B$. Consider an
ultraproduct $\left(  X_{\alpha}\right)  _{D}$. Since $X_{\alpha
}\hookrightarrow X$ for all $\alpha<\varkappa$, it may be assumed that
$\left(  X_{\alpha}\right)  _{D}\hookrightarrow\left(  X\right)  _{D}$ and,
hence, $\dim\left(  \left(  X_{\alpha}\right)  _{D}\right)  \leq\dim\left(
\left(  X\right)  _{D}\right)  $. So, to prove the theorem it is enough to
show that $\dim\left(  \left(  X_{\alpha}\right)  _{D}\right)  \geq
\varkappa^{\tau}$.

From $\tau$-regularity of $D$ it follows that there exists a subset $G\subset
D$, $\operatorname{card}(G)=\tau$, such that any $i\in I$ belongs only to a
finite number of sets $e\in D$.

Define on $G$ a linear order, say, $\ll$.

Let $g:G\rightarrow A$ be a function.

Let a function $f_{g}:I\rightarrow B$ be given by
\[
f_{g}\left(  i\right)  =\operatorname{span}\{a_{g\left(  e_{k}\left(
i\right)  \right)  }:k<n_{i}\}=B_{i},
\]
where $n_{i}$ is a cardinality of a set of those sets $e\in G$, to which $i$
belongs;
\[
\{e_{0}\left(  i\right)  \ll...\ll e_{n_{i}}\left(  i\right)  \}
\]
be a list of them.

Let $h:G\rightarrow A$ be any other function: $h\neq g$. Then there exists
$e\in G$ such that $h\left(  e\right)  \neq g\left(  e\right)  $.
Corresponding functions $f_{g}$ and $f_{h}$ are differ also. Indeed, for any
$i\in e$ the set $e$ is contained in a finite sequence $e_{0}\left(  i\right)
\ll...\ll e_{n_{i}}\left(  i\right)  $ of all sets that contain $i$. If its
order number in a such sequence is equal to $k$, then
\[
f_{g}\left(  i\right)  =\operatorname{span}\{...,g\left(  e_{k}\right)
,...\}\neq\operatorname{span}\{...,h\left(  e_{k}\right)  ,...\}=f_{h}\left(
i\right)  .
\]

Note that $e\in D$ and that $f_{g}\left(  i\right)  \neq f_{h}\left(
i\right)  $ for all $i\in e$.

Hence $f_{g}\left(  i\right)  $ and $f_{h}\left(  i\right)  $ generate
different elements $\frak{x}\left(  f\right)  $ and $\frak{x}\left(  h\right)
$ of $\left(  X_{\alpha}\right)  _{D}$; moreover $\left\|  \frak{x}\left(
f\right)  -\frak{x}\left(  h\right)  \right\|  \geq\varepsilon$.

Since $\operatorname{card}(^{G}A)=\varkappa^{\tau}$, different elements
$f,h\ $of the set $^{G}A$ of all functions from $G$ to $A$ generate different
elements $\frak{x}\left(  f\right)  $ and $\frak{x}\left(  h\right)  $ of the
ultraproduct $\left(  X_{\alpha}\right)  _{D}$ and, in addition, $\left\|
\frak{x}\left(  f\right)  -\frak{x}\left(  h\right)  \right\|  \geq
\varepsilon$, it is clear that $\dim\left(  \left(  X_{\alpha}\right)
_{D}\right)  \geq\varkappa^{\tau}$.
\end{proof}

\begin{remark}
It may be proved that for any infinite-dimensional Banach space $X$ and any
countably incomplete ultrafilter $D$,
\[
\dim((X)_{D})=(\dim((X)_{D}))^{\omega}.
\]
Hence, by using ultrapowers, cannot be obtained any space $(X)_{D}$, which
dimension $\varkappa$ is of countable confinality (i.e. such that
$\operatorname{cf}(\varkappa)=\omega$; e.g., $\omega_{\omega}$; $\omega
_{\omega_{\omega}}$, $\gimel_{\omega}$ and so on).
\end{remark}

\subsection{Elementary equivalent Banach spaces}

Let $X,Y\in\mathcal{B}$. These spaces are said to be \textit{elementary
equivalent} (in symbol: $X\equiv Y$) if there exists such ultrafilter $D$ that
ultrapowers $(X)_{D}$ and $(Y)_{D}$ are isometric.

It easy to prove that the relation $\equiv$ is reflexive and symmetric:%
\[
X\equiv X\text{ \ and \ }X\equiv Y\Leftrightarrow Y\equiv X.
\]
To show its transitivity it needs much more work.

Namely, it will be used the result that was presented without proof in
\cite{HeiMam}, which was called there \textit{the Keisler-Shelah theorem}. In
\cite{HeiMam} this result was formulated as follows:

\begin{itemize}
\item \textit{For every Banach space }$X$\textit{ and a pair of ultrafilters
}$E$\textit{ and }$D$\textit{ there exists a such ultrafilter }$G$\textit{
that ultrapowers }$((X)_{D})_{G}$\textit{ and }$((X)_{E})_{G}$\textit{ are isometric.}
\end{itemize}

We shall need the more powerful result, which is a consequence of the S.
Shelah's one (cf. \cite{Sh}), and which will be called the \textit{full
Keisler-Shelah theorem.}

\begin{theorem}
(\textbf{Keisler-Shelah}) Let $\varkappa$ be a cardinal. Let $X$ be a\textit{
Banach space and a pair of ultrafilters }$E$\textit{ and }$D$\textit{ be such
that }$\max\dim\{(X)_{D},(X)_{E}\}<\tau$, where $\tau$ is the least cardinal
such that $\varkappa^{\tau}>\varkappa$. Then \textit{there exists an
ultrafilter }$G$\textit{ over the set of cardinality }$\varkappa$, which does
not depend on $X$, $E$ and $D$, such \textit{that ultrapowers }$((X)_{D})_{G}%
$\textit{ and }$((X)_{E})_{G}$\textit{ are isometric.In other words, }%
$(X)_{D}\equiv(X)_{E}$ (and, in particular, $X\equiv\left(  X\right)  _{D}$
for every ultrafilter $D$)..
\end{theorem}

As a corollary we obtain the following result.

\begin{theorem}
The relation $X\equiv Y$ is an equivalence relation on the class $\mathcal{B}$
of all Banach spaces.
\end{theorem}

\begin{proof}
As it was noted before, the relation $\equiv$ is reflexive ($X\equiv X$) and
symmetric ($X\equiv Y\Leftrightarrow Y\equiv X$).

Assume that $X\equiv Y$ and $Y\equiv Z$. This means that there exists such
ultrafilters $D$ and $E$ that $\left(  X\right)  _{D}=\left(  Y\right)  _{D}$
and $\left(  Y\right)  _{E}=\left(  Z\right)  _{E}$ (the symbol $=$ means the
isometry). By the preceding theorem there exists such ultrafilter $G$ over the
set $\varkappa$, where $\varkappa^{\tau}>\tau$ and $\max\dim\{(Y)_{D}%
,(Y)_{E}\}<\tau$, that $((Y)_{D})_{G}=((Y)_{E})_{G}$. Hence, $((X)_{D}%
)_{G}=((Z)_{E})_{G}$.

Since $X\equiv((X)_{D})_{G}$, there exists such ultrafilter $F$ that $\left(
X\right)  _{F}=(((X)_{D})_{G})_{F}$. Notice that $F$ depends only on dimension.

So, one may chose $F$ such that $\left(  Z\right)  _{F}=(((Z)_{E})_{G})_{F}$.

Surely, this implies that $\left(  X\right)  _{F}=\left(  Z\right)  _{F}$,
i.e., that $X\equiv Z$.
\end{proof}

Thus, every Banach space $X$ generates a class (of elementary equivalence)%
\[
X^{\xi}=\{Y\in\mathcal{B}:X\equiv Y\}.
\]

The second result, presented in \cite{HeiMam} (also, without proof) is the
result that was called in \cite{HeiMam} \textit{the L\"{o}wenheim-Skolem
theorem. }This result for its formulation needs the following

\begin{definition}
Let $X$ be a Banach space; $Y$ be its subspace ($Y\hookrightarrow X$).

$Y$ is said to be an elementary subspace of $X$ ( in symbols: $Y\prec X$) if
there exists an ultrafilter $D$ and an isometry $i:(Y)_{D}\rightarrow\left(
X\right)  _{D}$ such that $i\circ d_{Y}=d_{X}\circ id_{Y}$, where
$id_{Y}:Y\hookrightarrow X$ is the identical embedding; $d_{Y}$ and $d_{X}$
are canonical embeddings of $Y$ and $X$ respectively into their ultrapowers.
\end{definition}

In other words, $Y\prec X$ if for an ultrafilter $D$ and an isometry
$i:(Y)_{D}\rightarrow\left(  X\right)  _{D}$ the diagram%
\[%
\begin{array}
[c]{ccc}%
(Y)_{D} & \overset{i}{\rightarrow} & \text{ \ }(X)_{D}\\
\uparrow &  & \uparrow\\
Y & \hookrightarrow &  X
\end{array}
,
\]
where vertical arrows denote respective canonical embeddings commutes\textit{.}

After this definition one can formulate the L\"{o}wenheim-Skolem theorem
\cite{HeiMam}:

\begin{theorem}
(\textbf{L\"{o}wenheim-Skolem}) \textit{For every Banach space }$X$\textit{ of
infinite dimension }$\varkappa$\textit{ and each subset }$A\subset X$\textit{
of cardinality }$\operatorname*{card}A<\varkappa$\textit{ there exists an
elementary subspace }$Y_{A}\prec X$\textit{ of dimension }$\dim Y_{A}%
=\max\{\omega,\operatorname*{card}A\}$\textit{ such that }$A\subset Y_{A}$\textit{.}
\end{theorem}

Certainly, if $Y\prec X$ then $Y\equiv X$. So, the following result is true.

\begin{theorem}
For every infinite-dimensional Banach space $X$ the corresponding class
$X^{\xi}$ contains spaces of arbitrary given dimension.

If $X$ is of finite dimension, $X^{\xi}$ contains exactly one member -- the
space $X$ itself.
\end{theorem}

\begin{proof}
From previous results it follows the second part of the theorem. If $X$ is of
infinite dimension (say, $\varkappa$), then, for every infinite cardinal
$\tau$ and regular ultrafilter $D$ over $\tau$, the ultrapower $\left(
X\right)  _{D}$ belongs to $X^{\xi}$. Its dimension $\dim$ $\left(  X\right)
_{D}=\varkappa^{\tau}$ may be maiden arbitrary large by choosing of $D$..

By the L\"{o}wenheim-Skolem theorem, for every infinite cardinal
$\kappa<\varkappa^{\tau}$ there exists an elementary subspace $X_{0}\prec$
$\left(  X\right)  _{D}$ of dimension $\kappa$.

Since $\kappa$ is arbitrary and $X_{0}\equiv X$, the theorem is done.
\end{proof}

So, any such class $X^{\xi}$ contains at least one separable space and it may
be considered the set
\[
\xi(\mathcal{B})=\{X^{\xi}\cap\mathcal{B}_{\omega}:X\in\mathcal{B}\}.
\]

Since there is one-to-one correspondence between classes $X^{\xi}$ and sets
$X^{\xi}\cap\mathcal{B}_{\omega}$, it will be assumed (for convenience) that
elements (members) of $\xi(\mathcal{B})$ are just classes $X^{\xi}$.

As it was already noted, both theorems -- the Keisler-Shelah and the
L\"{o}wenheim-Skolem were appeared in literature in the Banach space setting
without proofs. For the completeness we present in the following subsection
their proofs, which use some model theory

\subsection{Some model theory}

To express the notion of a Banach space in the logic consider a
\textit{language }$\frak{L}$ - a set of symbols, which includes besides of
logical symbols (connectives $\&,\vee,\urcorner$, quantifiers $\forall
,\exists$ and variables $u,v,x_{1},x_{2},\ $etc) non-logical primary symbols
as well:

\begin{itemize}
\item  A binary functional symbol $+$;

\item  A countable number of unary functional symbols $\{f_{q}:q\in
\mathbb{Q}\}$;

\item  An unary predicate symbol $B$.
\end{itemize}

Any Banach space $\frak{X}=\left\langle X,\left\|  \cdot\right\|
\right\rangle $ may be regarded as a model for all logical propositions (or
formulae) that are satisfied in $X$. A set $\left|  \frak{X}\right|  =X$ is
called a support (or an absolute) of a model $\frak{X}$, in which all non
logical symbols of the language $\frak{L}$ are interpreted as follows:

\begin{itemize}
\item $+^{\frak{X}}$ is interpreted as the addition of vectors from $X$;

\item ($f_{q})^{\frak{X}}$ for a given $q\in\mathbb{Q}$ is interpreted as a
multiplication of vectors from $X$ by a rational scalar $q$

\item $B^{\frak{X}}$ is interpreted as the unit ball $B(X)=\{x\in X:\left\|
x\right\|  \leq1\}$.
\end{itemize}

The language $\frak{L}$ was introduced by J. Stern \cite{stern}, who used it
for examine the Banach space\ theory in the \textit{first order logic}. In
this logic any quantifier acts only on elements of a support of corresponding
structure (e.g., a quantification over either all countable subsets of $X$ or
over all natural $n$'s is forbidden); only formulae that contains only finite
strings of symbols are admissible.

For a Banach space $X$ (which is identified with the model $\frak{X}$) define
a set $\operatorname*{Th}(X)$ of all formulae of the first order logic,
expressed in the language $\frak{L}$, which are satisfied in $\frak{X}$ (in
this case we shall write $\frak{X}\vDash\operatorname*{Th}(X)$). Certainly, in
a general case $\frak{X}$ is not the unique model of $\operatorname*{Th}(X)$:
there may be exist other structures, different from $\frak{X}$ (say,
$\frak{A}$) such that $\frak{A}$ satisfies the same formulae as $\frak{X}$ (in
symbol: $\frak{A}\vDash\operatorname*{Th}(X)$). Of course, $\frak{A}$ and
$\frak{X}$ cannot be distinguished by tools of the first order logic.

There is a method to produce such models $\frak{A}$. It is called '' the
model-theoretical ultrapower'' and is denoted by $\Pi_{D}\frak{X}$, where $D$
is an ultrafilter over a set $I$ and $\Pi_{D}\frak{X}$ is a model (or,
equivalently, structure) for $\frak{L}$, whose absolute (or support) is the
set-theoretical ultrapower of set $X$ (i.e. a quotient $\Pi_{D}X$ of the
Carthesian product $\Pi_{I}X$ of $\operatorname*{card}I$ copies of $X$ by the
equivalence relation: $\left(  x_{i}\right)  =_{D}\left(  y_{i}\right)  $ that
means that $\{i\in I:x_{i}=y_{i}\}\in D$), and non-logical primary symbols of
$\frak{L}$ are interpreted in $\Pi_{D}\frak{X}$ in the following way:

\begin{itemize}
\item  Addition of vectors and multiplying by a rational scalar are defined coordinate-wise;

\item $B^{\Pi_{D}\frak{X}}$ is interpreted by
\[
B^{\Pi_{D}\frak{X}}=\{\frak{x}=\left(  x_{i}\right)  _{D}\in\Pi_{D}%
\frak{X}:\{i\in I:\left\|  x_{i}\right\|  \leq1\}\in D\}.
\]
\end{itemize}

It is well-known that $\Pi_{D}\frak{X}\vDash\operatorname*{Th}\left(
X\right)  $. However, $\Pi_{D}\frak{X}$ is not necessary a Banach space
because of it may contain nonstandard elements (non-zero elements with zero
norm or elements for which their norms, calculated by the Minkowski's
functional, are infinite).

J. Stern \cite{stern} suggested a procedure of elimination of nonstandard
elements, which looks as follows. Let $\frak{A}\vDash\operatorname*{Th}\left(
X\right)  $ be a (may be - nonstandard) model. Put%
\begin{align*}
\operatorname*{Fin}\left(  \frak{A}\right)   &  =\{\frak{a}\in\frak{A}%
:B^{\frak{A}}(n^{-1}\frak{a})\text{ \ for some }n\in\mathbb{N}\};\\
\operatorname*{Null}\left(  \frak{A}\right)   &  =\{\frak{a}\in\frak{A}%
:B^{\frak{A}}(n^{-1}\frak{a})\text{ \ for all }n\in\mathbb{N}\};\\
\left\|  \frak{a}\right\|  ^{\frak{A}} &  =\inf\{q\in\mathbb{Q}:B^{\frak{A}%
}(q^{-1}\frak{a})\}
\end{align*}
and define the procedure $\left[  \cdot\right]  :\frak{A}\rightarrow\left[
\frak{A}\right]  $ by%
\[
\left[  \frak{A}\right]  =(\operatorname*{Fin}\left(  \frak{A}\right)
/\operatorname*{Null}\left(  \frak{A}\right)  )\symbol{94},
\]
where $\symbol{94}$ means the completition by the norm $\left\|
\cdot\right\|  ^{\frak{A}}$.

By \cite{stern}, this procedure sends the model-theoretical ultrapower
$\Pi_{D}\frak{X}$ of a Banach space $X$, which is regarded as an $\frak{L}$
-model $\frak{X}$ to its Banach-space ultrapower $\left(  X\right)  _{D}.$

After the preparation we are ready to present proofs of both results mentioned before.

\begin{proof}
\textit{The Keisler-Shelah theorem}\textbf{. }

Let\textbf{ }$X$ be a Banach space regarding as an $\frak{L}$-model $\frak{X}%
$. Let $D$, $E$ be ultrafilters. Consider model-theoretical ultrapowers
$\Pi_{D}\frak{X}$ and $\Pi_{E}\frak{X}$. Since, as it was noted before,
$\frak{X}\vDash\operatorname*{Th}\left(  X\right)  $; $\Pi_{D}\frak{X}%
\vDash\operatorname*{Th}\left(  X\right)  $ and $\Pi_{E}\frak{X}%
\vDash\operatorname*{Th}\left(  X\right)  $ as well, by the result of
J.Keisler \cite{kei} (in assumption the GCH) and of S. Shelah \cite{Sh} (in a
general case) there exists a such ultrafilter $F$ that $\Pi_{F}(\Pi
_{D}\frak{X})=\Pi_{F}(\Pi_{E}\frak{X})$. Moreover, $F$ may be chosen in a such
way, that it depends only of maximal cardinality $\tau=\max\dim\{\Pi
_{E}\frak{X},\Pi_{D}\frak{X}\}$ (to be an ultrafilter over a set $\varkappa
$where $\varkappa$ is the minimal cardinal such that $\tau^{\varkappa
}>\varkappa$ ; see \cite{Sh}.

Hence
\[
\left[  \cdot\right]  :\Pi_{F}(\Pi_{D}\frak{X})\rightarrow\left[  \Pi_{F}%
(\Pi_{D}\frak{X})\right]  =(\left[  \Pi_{D}\frak{X}\right]  )_{F}=\left(
\left(  X\right)  _{D}\right)  _{F}.
\]

From the other hand,%
\[
\left[  \cdot\right]  :\Pi_{F}(\Pi_{E}\frak{X})\rightarrow\left[  \Pi_{F}%
(\Pi_{E}\frak{X})\right]  =(\left[  \Pi_{E}\frak{X}\right]  )_{F}=\left(
\left(  X\right)  _{E}\right)  _{F}.
\]

So, $\left(  \left(  X\right)  _{D}\right)  _{F}=\left(  \left(  X\right)
_{E}\right)  _{F}$ and the theorem is proved.
\end{proof}

\begin{proof}
\textit{The L\"{o}wenheim-Skolem theorem. }

For $\varkappa=\omega$ the result is obvious: as $Y_{A}$ it may be chosen the
space $X$ itself.

Let $\varkappa>\omega$. Let $X$ be a Banach space of dimension $\varkappa
>\omega$; $A\subset X$ be a subset of cardinality $\tau$; $\omega\leq
\tau<\varkappa$. Consider $X$ as an $\frak{L}$-model $\frak{X}$. By the (model
theoretical) L\"{o}wenheim-Skolem theorem there exists an submodel (say
$\frak{Y}$) of $\frak{X}$ such that its absolute $\left|  \frak{Y}\right|  $
contains $A$, $\frak{Y}\vDash\operatorname*{Th}\left(  X\right)  $ and
$\operatorname*{card}\left|  \frak{Y}\right|  =\operatorname*{card}A$.
Moreover, $\frak{Y}$ is an elementary submodel of $\frak{X}$, i.e. there are
an ultrafilter $D$ and a model-theoretical isomorphism $m:\Pi_{D}%
\frak{Y}\rightarrow\Pi_{D}\frak{X}$ such that the following diagram commutes:%
\[%
\begin{array}
[c]{ccc}%
\Pi_{D}\frak{Y} & \overset{m}{\rightarrow} & \Pi_{D}\frak{X}\\
\uparrow &  & \uparrow\\
\frak{Y} & \hookrightarrow & \frak{X}%
\end{array}
\]

Using the Stern's procedure $\left[  \cdot\right]  $ it immediately follows
that $\left[  \frak{Y}\right]  \prec\left[  \frak{X}\right]  $ and that
$\dim\left[  \frak{Y}\right]  =\dim A.$

Put $Y_{A}=\left[  \frak{Y}\right]  $. Clearly, $A\subset Y_{A}\prec X$.
\end{proof}

\subsection{Topological properties of $\xi(\mathcal{B})$}

In this subsection it will be used another definition of the notion of
elementary equivalence, that belongs to S. Heinrich and C.W.Henson, see
\cite{hh}. Namely,

\begin{itemize}
\item \textit{Banach spaces }$X$\textit{ and }$Y$\textit{ are elementary
equivalent, }$X\equiv Y$\textit{, if and only if for any natural }$n$\textit{,
any }$\varepsilon>0$\textit{ and any chain of finite dimensional subspaces
}$A_{1}\hookrightarrow A_{2}\hookrightarrow...\hookrightarrow A_{n}$\textit{
of }$X$\textit{ there exists a chain of finite dimensional subspaces }%
$B_{1}\hookrightarrow B_{2}\hookrightarrow...\hookrightarrow B_{n}$\textit{ of
}$Y$\textit{ such that there exists an isomorphism }$u:B_{n}\rightarrow A_{n}%
$\textit{ with }$u(B_{i})=A_{i}$\textit{ for all }$i<n$\textit{, which
satisfies }$\left\|  u\right\|  \left\|  u^{-1}\right\|  <1+\varepsilon
$\textit{, corresponding restrictions }$u\mid_{B_{i}}$\textit{ have the same
estimates of norms and for any finite dimensional subspace }$A_{n+1}$\textit{,
}$A_{n}\hookrightarrow A_{n+1}\hookrightarrow X$\textit{ there exists a
subspace }$B_{n+1}$\textit{, }$B_{n}\hookrightarrow B_{n+1}\hookrightarrow
Y$\textit{ with an isomorphism }$\overline{u}$\textit{ between }$B_{n+1}%
$\textit{ and }$A_{n+1}$\textit{, which restriction to }$B_{n}$\textit{ is
equal to }$u$\textit{ and which satisfies the estimate }$\left\|  u\right\|
\left\|  u^{-1}\right\|  <1+\varepsilon$\textit{ and, conversely, for any
B}$_{n+1}$\textit{, }$B_{n}\hookrightarrow B_{n+1}\hookrightarrow Y$\textit{
there exists }$A_{n+1}$\textit{, }$A_{n}\hookrightarrow A_{n+1}\hookrightarrow
X$\textit{ and an isomorphism }$\overline{u}$\textit{ between B}$_{n+1}%
$\textit{ and A}$_{n+1}$\textit{ which restriction to }$B_{n}$\textit{ is
equal to }$u$\textit{ and which norm satisfies the same estimate.}
\end{itemize}

Let $\frak{A}\subset\xi(\mathcal{B}).$ Recall that we assume elements of
$\frak{A}$ are classes $X^{\xi}$.

The set $\frak{A}$ is said to be \textit{ultraclosed} if for any family
$\{X_{i}:i\in I\}$, where $(X_{i})^{\xi}$ belongs to $\frak{A}$, their
ultraproduct $Z=(X_{i})_{D}$ generates a class $Z^{\xi}\in\frak{A}$.

It will be said that $\frak{A}$ is \textit{double ultraclosed} if both sets
$\frak{A}$ and $\xi(\mathcal{B})\backslash\frak{A}$ are ultraclosed.

Define on $\xi(\mathcal{B})$\ a topology $\mathcal{T}_{\xi}$ by choosing for
its base of open sets a set $\frak{W}$ of all double ultraclosed sets
$A\subset$.$\xi(\mathcal{B})$.

\begin{theorem}
The topological space $\left\langle \xi(\mathcal{B}),\mathcal{T}_{\xi
}\right\rangle $ is totally disconnected, Hausdorff and compact.
\end{theorem}

\begin{proof}
The first assertion about $\xi(\mathcal{B})$ follows from the definition of
the topology.

\textit{Compactness}. Let $\frak{B}\in\frak{A}$ has the property: for any
finite subset $\frak{B}_{0}$ of $\frak{B}$ its intersection $\cap\frak{B}_{0}$
is nonempty. It may be defined such ultrafilter $D$ that for any family
$(X_{i})_{i\in\frak{B}}$, where each $X_{i}$ generates a member of $\frak{B}$,
their ultraproduct $(X_{i})_{D}$ generates a member of the intersection $\cap
B$. Hence, $\xi(\mathcal{B})$ has the finite intersection property and,
consequently, is compact.

Let us show that $\left\langle \xi(\mathcal{B}),\mathcal{T}_{\xi}\right\rangle
$ is \textit{Hausdorff}.

Let $X^{\xi}$, $Y^{\xi}$ be distinct classes. There exist natural $n$, real
$\varepsilon>0$ and a chain $A_{1}\hookrightarrow A_{2}\hookrightarrow
..\hookrightarrow A_{n}$ of finite dimensional subspaces of $X$ such that for
any $Y\in Y^{\xi}$ and for any chain $B_{1}\hookrightarrow B_{2}%
\hookrightarrow...\hookrightarrow B_{n}$ of finite dimensional subspaces of
$Y$, which is $(1+\varepsilon)$-isomorphic to a chain $A\hookrightarrow
A\hookrightarrow...\hookrightarrow A_{n}$ there exists $A_{n+1}$,
$A_{n}\hookrightarrow A_{n+1}\hookrightarrow X$ such that $Y_{0}$ does not
contain any $B_{n+1}$, $B_{n}\hookrightarrow B_{n+1}\hookrightarrow Y$, that
extended $B_{1}\hookrightarrow B_{2}\hookrightarrow...\hookrightarrow B_{n}$
to the chain which is $(1+\varepsilon)$-isomorphic to $A_{1}\hookrightarrow
A_{2}\hookrightarrow...\hookrightarrow A_{n+1}$ (or respectively such
$B_{n+1}$, $B_{n}\hookrightarrow B_{n+1}\hookrightarrow Y$ that no one of
$A_{n+1}$, $A_{n}\hookrightarrow A_{n+1}\hookrightarrow X$ forms a desired extension).

Let $\frak{Y}_{k,n}$ be a class of all $Z\in\mathcal{B}$ which have the
property:

Any $(1+1/k)$-isomorphism between given chains $A_{1}\hookrightarrow
A_{2}\hookrightarrow...\hookrightarrow A_{n}$ of $X$ and $B_{1}\hookrightarrow
B_{2}\hookrightarrow...\hookrightarrow B_{n}$ of $Z$ cannot be extended to a
$(1+1/k)$-isomorphism between $A_{n+1}$, $A_{n}\hookrightarrow A_{n+1}%
\hookrightarrow X$ and $B_{n+1}$, $B_{n}\hookrightarrow B_{n+1}\hookrightarrow
Z$ and conversely.

It is clear that the class $\frak{G}=\cup\{\frak{Y}_{k,n}:n,k<\infty\}$\ is
doubly ultraclosed. Certainly, this class contains $Y^{\xi}$. The class
$X^{\xi}$ belongs to $\xi(\mathcal{B})\backslash\frak{G}$.
\end{proof}

\section{Classification of Banach spaces by finite equivalence}

\begin{definition}
Let $X$, $Y$ be Banach spaces. $X$ is said to be \textit{finitely
representable} in $Y$, shortly: $X<_{f}Y$, if for every $\varepsilon>0$ and
for every finite dimensional subspace $A$ of $X$ there is a subspace $B$ of
$Y$ and an isomorphism $u:A\rightarrow B$ such that $\left\|  u\right\|
\left\|  u^{-1}\right\|  <1+\varepsilon$. $X$ and $Y$ are said to
be\textit{\ finitely equivalent}, $X\sim_{f}Y$, if $X<_{f}Y$ and $Y<_{f}X$.
So, every Banach space $X$ generates a \textit{class of finite equivalence}
\[
X^{f}=\{Y\in\mathcal{B}:X\sim_{f}Y\}.
\]
\end{definition}

For any infinite dimensional Banach space $X$ the corresponding class $X^{f}$
is proper: it contains spaces of arbitrary large dimension.

\begin{remark}
It \ is well known that $X<_{f}Y$ if and only if $X$ is isometric to a
subspace of some ultrapower $\left(  Y\right)  _{D}$.
\end{remark}

\begin{theorem}
For any infinite dimensional Banach space $X$ the class $X^{f}$ contains at
least one space $Y_{\varkappa}$ of dimension $\varkappa$ for each infinite
cardinal number $\varkappa$ (particularly, at least one separable space
$Y_{\omega}$).In other words, for any cardinal $\varkappa\geq\omega$,
$X^{f}\cap\mathcal{B}_{\varkappa}\neq\varnothing$
\end{theorem}

\begin{proof}
Let $\dim(X)=\tau$. If $\varkappa=\omega$ then we choose a countable sequence
$\{A_{i}:i<\omega\}$ of finite dimensional subspaces of $X$ which is dense in
a set $H(X)$ of all different finite dimensional subspace of $X$ (isometric
subspaces in $H(X)$are identified), equipped with a metric topology, which is
induced by the Banach-Mazur distance. Clearly, $X_{0}=\operatorname{span}%
\{A_{i}:i<\omega\}\hookrightarrow X$ is separable and is finitely equivalent
to $X.$

If $\varkappa\leq\tau$ then as a representative of $X^{f}\cap\mathcal{B}%
_{\varkappa}$ may be chosen any subspace of $X$ of dimension $\tau
\ $that\ contains a subspace $X_{0}$.

If $\tau<\varkappa$ then, by the theorem 8 it may be found an ultrapower
$(X)_{D}$ of dimension $\geq\varkappa$. Certainly, $(X)_{D}\sim_{f}X$. Choose
any subspace of $(X)_{D}$ of dimension $\varkappa$, which contains a subspace
$d_{X}X$.
\end{proof}

\begin{remark}
This result is an obvious consequence of the \textit{L\"{o}wenheim-Skolem theorem.}
\end{remark}

For any two Banach spaces $X$, $Y$ their \textit{Banach-Mazur distance }is
given by%
\[
d(X,Y)=\inf\{\left\|  u\right\|  \left\|  u^{-1}\right\|  :u:X\rightarrow Y\},
\]
where $u$ runs all isomorphisms between $X$ and $Y$ and is assumed, as usual,
that $\inf\varnothing=\infty$.

It is well known that $\log d(X,Y)$ forms a metric on each class of isomorphic
Banach spaces, where almost isometric Banach spaces are identified.

Recall that Banach spaces $X$ and $Y$ are \textit{almost isometric} if
$d(X,Y)=1$. Surely, any almost isometric finite dimensional Banach spaces are isometric.

The set $\frak{M}_{n}$ of all $n$-dimensional Banach spaces, equipped with
this metric, is the compact metric space, called \textit{the Minkowski
compact} $\frak{M}_{n}$.

The disjoint union $\cup\{\frak{M}_{n}:n<\infty\}=\frak{M}$ is a separable
metric space, which is called the \textit{Minkowski space}.

Consider a Banach space $X$. Let $H\left(  X\right)  $ be a set of all its
different finite dimensional subspaces (isometric finite dimensional subspaces
of $X$ in $H\left(  X\right)  $ are identified). Thus, $H\left(  X\right)  $
may be regarded as a subset of $\frak{M}$, equipped with a restriction of the
metric topology of $\frak{M}$.

Of course, $H\left(  X\right)  $ need not to be a closed subset of $\frak{M}$.
Its closure in $\frak{M}$ will be denoted by $\overline{H\left(  X\right)  }$.
From definitions it follows that $X<_{f}Y$ if and only if $\overline{H\left(
X\right)  }\subseteq\overline{H\left(  Y\right)  }$ and therefore, $X\sim
_{f}Y$ if and only if $\overline{H\left(  X\right)  }=\overline{H\left(
Y\right)  }$. Thus, there is a one to one correspondence between classes of
finite equivalence
\[
X^{f}=\{Y\in\mathcal{B}:X\sim_{f}Y\}
\]
and closed subsets of $\frak{M}$ of kind $\overline{H\left(  X\right)  }$.
Indeed, all spaces $Y$ from $X^{f}$ have the same set $\overline{H\left(
X\right)  }$. This set, uniquely determined by $X$ (or, equivalently, by
$X^{f}$), will be denoted $\frak{M}(X^{f})$ and will be referred to as
\textit{the Minkowski's base of the class} $X^{f}$.

Let $f(\mathcal{B})$ be a set of all different classes of kind $X^{f}$.

To avoid set-theoretical difficulties note that each class $X^{f}$ is in
one-to-one correspondence with the set $\frak{M}(X^{f})$ and $f(\mathcal{B})$
may be regarded as the set of sets of kind $\frak{M}(X^{f})$.

The following result is obvious.

\begin{theorem}
$\operatorname*{card}f(\mathcal{B})=\frak{c}$ $\left(  =\exp\left(
\omega\right)  \right)  $.

\begin{proof}
There is just the continuum number of closed subsets of the Minkowski's space
$\frak{M}$.

At the same time, there is the continuum number of pairwice distinct classes
of finite equivalence, e.g., classes $(l_{p})^{f}$ ($1\leq p\leq\infty$) .
\end{proof}
\end{theorem}

\subsection{Topological properties of $\frak{M}(X^{f})$.}

Let $X$ be a finite dimensional Banach space; $\dim X=n$. In this subsection
it will be shown that for $n\geq2$ either $\operatorname{card}\frak{M}%
(X^{f})=\frak{c}$ or $X$ is isomorphic to the Euclidean space $l_{2}^{(n)}$
(in this case $\operatorname{card}\frak{M}(X^{f})=n+1$).

Of course, $\operatorname{card}\frak{M}(X^{f})=3$ for any 2-dimensional space
$X$ since each of its subspaces is either $X$ itself or is $0$- or $1$-dimensional.

Along with $H(X)$ it will be considered a set $G(X)$ of all distinct
finite-dimensional subspaces of $X$ ($A,B\in G(X)$ are considered as distinct
if either they are of different dimension or $A\cap B\neq A$). A subset%
\[
G_{m}(X)=\{A\in G(X):\dim A=m\},
\]
endowed with the metric%
\[
\Theta(A,B)=\max\{r(A,B),r(B,A)\},
\]
where%
\[
r(A,B)=\sup\{\inf\{\left\|  a-b\right\|  :\left\|  b\right\|  =1\}:\left\|
a\right\|  =1\}
\]
is called \textit{the Grassman manifold}. It is known that $\left\langle
G_{m}(X),\Theta\right\rangle $ is a compact connected metric space for any
finite dimensional $X\in\mathcal{B}$ and every $m<\dim(X)$.

\begin{theorem}
$\left\langle H_{m}(X),\varrho\right\rangle $ is a connected (compact, metric) space.
\end{theorem}

\begin{proof}
Let $(A_{i})_{i<\infty}\subset G_{m}(X)$ be a fundamental sequence with
respect to the metric $\Theta$, i.e. such that $\lim\nolimits_{i,k\rightarrow
\infty}\Theta(A_{i},A_{k})=0$.

Because of compactness of $G_{m}(X)$ there is such $A\in G_{m}(X)$ and a
subsequence $(A_{i_{k}})_{k<\infty}\subset(A_{i})_{i<\infty}$ that
$\Theta(A_{i_{k}},A)\leq m^{-(k+1)}=\left(  \dim A\right)  ^{-k-1}$.

Let $(e_{j})_{j<m}$ be the \textit{Auerbach basis} of $A$, i.e., such that for
any element $a\in A$, $a=\sum\nolimits_{j=1}^{m}a_{j}e_{j}$, the following
inequality holds:%
\[
\max\nolimits_{1\leq j\leq m}\left|  a_{j}\right|  \leq\left\|  a\right\|
\leq\sum\nolimits_{1\leq j\leq m}\left|  a_{j}\right|  .
\]
Chose in $A_{i_{k}}$ vectors $(e_{j}^{k})$ of norm one such that for all
$j=1,2,...,m$
\[
\left\|  e_{j}-e_{j}^{k}\right\|  \leq\Theta(A_{i_{k}},A)\leq m^{-(k+1)}%
\]
Clearly,
\[
\sum\nolimits_{j=1}^{m}\left\|  e_{j}-e_{j}^{k}\right\|  \leq mm^{-(k+1)}%
=m^{-k}<1.
\]
Hence, by the Krein-Milman-Rutman theorem, $(e_{j}^{k})_{j=1}^{m}$ is a basis
of $A_{i_{k}}$.

Let an isomorphism $u_{k}:A\rightarrow A_{k}$ ($k<\infty$) be given by
$u_{k}\left(  e_{j}\right)  =e_{j}^{k}$ ($j=1,2,...,m$).

Since $(e_{j})_{j<m}$ is the Auerbach basis of $A$,
\begin{align*}
\left\|  u_{k}\right\|   &  =\sup\left\|  \sum\nolimits_{j=1}^{m}a_{j}%
e_{j}\right\|  /\left\|  \sum\nolimits_{j=1}^{m}a_{j}e_{j}^{k}\right\|  \\
&  =\sup\left\|  \sum\nolimits_{j=1}^{m}a_{j}e_{j}\right\|  /\left\|
\sum\nolimits_{j=1}^{m}a_{j}e_{j}+\sum\nolimits_{j=1}^{m}a_{j}\left(
e_{j}^{k}-e_{j}\right)  \right\|  \\
&  \leq\sup\left\|  \sum\nolimits_{j=1}^{m}a_{j}e_{j}\right\|  /\left(
\left\|  \sum\nolimits_{j=1}^{m}a_{j}e_{j}\right\|  -\left\|  \sum
\nolimits_{j=1}^{m}a_{j}\left(  e_{j}^{k}-e_{j}\right)  \right\|  \right)  \\
&  \leq\sup\left(  1-\sum\nolimits_{j=1}^{m}\left|  a_{j}\right|  \left\|
e_{j}^{k}-e_{j}\right\|  /\left\|  \sum\nolimits_{j=1}^{m}a_{j}e_{j}\right\|
\right)  ^{-1}\\
&  \leq\sup\left(  1-m^{-k-1}\sum\nolimits_{j=1}^{m}\left|  a_{j}\right|
/\max\nolimits_{1\leq j\leq m}\left|  a_{j}\right|  \right)  ^{-1}\leq\left(
1-m^{-k}\right)  ^{-1}.
\end{align*}

Similarly the norm of the inverse operator $u_{k}^{-1}$ is estimated:%
\begin{align*}
\left\|  u_{k}^{-1}\right\|   &  =\sup\left\|  \sum\nolimits_{j=1}^{m}%
a_{j}e_{j}^{k}\right\|  /\left\|  \sum\nolimits_{j=1}^{m}a_{j}e_{j}\right\|
\\
&  \leq\sup\left(  \left\|  \sum\nolimits_{j=1}^{m}a_{j}e_{j}\right\|
+\left\|  \sum\nolimits_{j=1}^{m}a_{j}\left(  e_{j}^{k}-e_{j}\right)
\right\|  \right)  /\left\|  \sum\nolimits_{j=1}^{m}a_{j}e_{j}\right\|  \\
&  \leq\sup\left(  \left\|  \sum\nolimits_{j=1}^{m}a_{j}e_{j}\right\|
+\sum\nolimits_{j=1}^{m}\left|  a_{j}\right|  m^{-k-1}\right)  /\left\|
\sum\nolimits_{j=1}^{m}a_{j}e_{j}\right\|  \\
&  \leq1+m^{-k}.
\end{align*}

Hence, $\left\|  u_{k}\right\|  \left\|  u_{k}^{-1}\right\|  \leq\left(
1+m^{-k}\right)  /\left(  1-m^{-k}\right)  \rightarrow1$ as $k\rightarrow
\infty$.

Hence, the convergence of sequences with respect to the metric $\Theta$
implies the convergence with respect to the metric $\varrho$. Moreover, this
implies that $H_{m}(X)$ is the connected space because it is a continuous
image of the connected space $G_{m}\left(  X\right)  $ (under the mapping that
'gluing' isometric subspaces).
\end{proof}

\begin{corollary}
If $X$ is not equal to $l_{2}^{\left(  n\right)  }$ and $\dim\left(  X\right)
\geq3$ then $\operatorname{card}H(X)=\frak{c}$.
\end{corollary}

\begin{proof}
Let $m<n$. Than either $\operatorname{card}H_{m}(X)=1$ or $\operatorname{card}%
H_{m}(X)=\frak{c}$. In the first case the condition dim$\left(  X\right)
\geq3$ yields that there exists an even $m<n$. According to the M. Gromov's
theorem \cite{Gr}, conditions $\operatorname{card}H_{m}(X)=1$, $m<n$ and $m$
is even imply that $X$ is isometric to the Euclidean space.
\end{proof}

\begin{corollary}
For every infinite dimensional Banach space $X$ and each $n<\omega$ the set
$H(X)$ is a connected subset of $\frak{M}_{n}$.
\end{corollary}

\begin{proof}
Fix $n<\omega$ and assume that $A$ and $B$ belong to different components of
$H_{n}(X)$. These spaces may be regarded as subspaces of $X$. Let
$C=\operatorname{span}\{A,B\}$. Certainly, $C\in H(X)$; $A\in H_{n}(C)$ and
$B\in H_{n}(C)$. However, $H_{n}(C)$ is a connected set according to the
theorem. Surely, this contradicts to our assumption.
\end{proof}

\begin{corollary}
For every infinite dimensional Banach space $X$ and each $n<\omega$ the set
$\frak{M}_{n}(X^{f})$ is a connected compact perfect space.
\end{corollary}

\begin{proof}
It is enough to notice that for every (non trivial) ultrafilter $D$ and every
$X\in\mathcal{B}$ its ultrapower $\left(  X\right)  _{D}$ satisfies $H\left(
\left(  X\right)  _{D}\right)  =\frak{M}(X^{f})$.

Hence, $H_{n}\left(  \left(  X\right)  _{D}\right)  =\frak{M}_{n}(X^{f})$ for
all $n<\omega$. So, by the previous corollary, $\frak{M}_{n}(X^{f})$ is
connected. Since $\frak{M}_{n}(X^{f})$ is compact it is perfect.
\end{proof}

\subsection{A synthesis of Banach spaces}

A problem of synthesis of Banach spaces is to describe sets $\frak{N}%
\in\frak{M}$ with the property: there is a Banach space $X_{\frak{N}}$ with
$\frak{M}(X^{f})=\frak{N}$ and to give a procedure for constructing such
$X_{\frak{N}}$. Obviously, necessary conditions on $\frak{N}$ that yields the
existence of $X_{\frak{N}}$ with $\frak{M}(X^{f})=\frak{N}$ are:

(\textbf{H}) If $A\in\frak{N}$ and $B\in H(A)$ then $B\in\frak{N}$.

(\textbf{A}$_{0}$) If $A,B\in\frak{N}$ then there exists $C\in\frak{N}$ \ such
that $A\in H(C)$ and $B\in H(C)$.

(\textbf{C}) $\frak{N}$ is closed in $\frak{M}$.

\begin{theorem}
Let $\frak{N}\in\frak{M}$. There exists a Banach space $X_{\frak{N}}$ with
$\frak{M}(X^{f})=\frak{N}$ if and only if $\frak{N}$ satisfies conditions
(\textbf{H}), (\textbf{A}$_{0}$), (\textbf{C}).
\end{theorem}

\begin{proof}
Assume that $\frak{N}\in\frak{M}$ is closed and satisfies conditions
(\textbf{H}) and (\textbf{A}$_{0}$). Define on $\frak{N}$ a partial order
assuming that $A<B$ if $A\in H\left(  B\right)  $ ($A$, $B\in\frak{N}$).

To every $A\in\frak{N}$\ corresponds a set
\[
A^{>}=\{B\in\frak{N}:A<B\}.
\]

From the property (\textbf{A}$_{0}$) it follows that $\mathcal{F}=\{A^{>}%
:A\in\frak{N}\}$ has the \textit{finite intersection property}:

\begin{itemize}
\item \textit{The intersection of any finite subset of }$\mathcal{F}$\textit{
is nonempty.}
\end{itemize}

Indeed, if $A_{0}$, $A_{1}$, ..., $A_{m-1}\in\frak{N}$ then there exists a
such $B_{m}\in\frak{N}$ \ that $A_{i}<B_{m}$ for all $i<m$. Hence, $\cap
_{i<m}A_{i}^{>}\supseteq B_{m}^{>}\neq\varnothing$. So, $\mathcal{F}$ may be
extended to some ultrafilter $D$ over $\frak{N}$.

Let $\frak{N}$ be considered as a set $\{A_{A}:A\in\frak{N}\}$, indexed by
itself.

By a standard technique of ultrapowers it follows that the ultraproduct
$W=(A_{A})_{D}$ of all spaces from $\frak{N}$ has desired properties.

Indeed, for every $B\in H(W)$ and every $\varepsilon>0$ there exists a space
$A\in\frak{N}$, which is $\left(  1+\varepsilon\right)  $-isomorphic to $B$.
Surely the converse is also true: $\frak{N}\subset H(W)$. Since $\frak{N}$ is
closed, $\frak{N}=H(W)$.

Conversely, every $\frak{N}\in\frak{M}$ that satisfy conditions (\textbf{H})
and (\textbf{A}$_{0}$) corresponds to some unique maximal centered system of
closed subsets of $\left\langle F(\mathcal{B}),\mathcal{T}\right\rangle $ and,
hence, to a class $(X_{\frak{N}})^{f}$. Any space of this class has desired properties.
\end{proof}

\subsection{Topological properties of $f(\mathcal{B})$}

Let $A\in\frak{M}(X^{f})$. Let $\mathcal{T}_{f}$ be a topology on
$f(\mathcal{B})$, which base of open sets is formed by finite unions and
intersections of sets of kind%
\[
U_{+}(A)=\{Y^{f}:A\in\frak{M}(X^{f})\};\text{ \ }U_{-}(A)=\{Y^{f}%
:A\notin\frak{M}(X^{f})\}.
\]

\begin{theorem}
A space $\left\langle f(\mathcal{B}),\mathcal{T}_{f}\right\rangle $ is a
totally disconnected Hausdorff compact topological space.
\end{theorem}

\begin{proof}
1. Since $U_{+}(A)=f(\mathcal{B})\backslash U_{-}(A)$, every open set
$U\in\mathcal{T}_{f}$ is also closed, i.e., $\left\langle f(\mathcal{B}%
),\mathcal{T}_{f}\right\rangle $ is totally disconnected.

2. If $X^{f}\neq Y^{f}$ then $\frak{M}(X^{f})\neq\frak{M}(Y^{f})$, i.e., there
exists such $A\in\frak{M}$ that $A\in\frak{M}(X^{f})$ and $A\notin
\frak{M}(Y^{f})$. From definitions it follows that $X^{f}\in U_{+}(A)$;
$Y^{f}\in U_{-}(A)$ and, hence, $\left\langle f(\mathcal{B}),\mathcal{T}%
_{f}\right\rangle $ is Hausdorff.

3. The compactness of $\left\langle f(\mathcal{B}),\mathcal{T}_{f}%
\right\rangle $ is a consequence of the theorem on synthesis. Indeed, let
$\mathcal{T}_{0}\subset\mathcal{T}_{f}$ be a centered system, i.e., every
finite collection of sets of $\mathcal{T}_{0}$ has the nonempty intersection.
Then $\mathcal{T}_{0}$ may be extended to an ultrafilter $D$ and the
ultrapower $(\mathcal{T}_{0})_{D}=W$ generates a class $W^{f}$.

Immediately, the intersection $\cap\mathcal{T}_{0}$ contains $W^{f}$ and,
hence, is nonempty as well. Since $\mathcal{T}_{0}$ is arbitrary,
$\left\langle f(\mathcal{B}),\mathcal{T}_{f}\right\rangle $ is compact.
\end{proof}

\begin{remark}
The same result may be obtained in other way: sets of kind $U_{+}(A)$ and
$U_{-}(A)$ generate a Boolean algebra (say, $\frak{B}$), which Stonian space
$\frak{S}\left(  \frak{B}\right)  $ is a totally disconnected Hausdorff
compact space according to the well-known M. Stone's theorem \cite{stone}.

Obviously, $\frak{S}\left(  \frak{B}\right)  $ and $\left\langle
f(\mathcal{B}),\mathcal{T}_{f}\right\rangle $ are homeomorphic.
\end{remark}

\section{Classification of Banach spaces by strong finite equivalence}

\begin{definition}
Let $X$, $Y\in\mathcal{B}$. $X$ is said to be \textit{strongly finitely
representable} in $Y$; shortly: $X<_{\phi}Y$, if $H(X)\subseteq H(Y)$. $X,Y$
$\in\mathcal{B}$ are said to be \textit{strongly finitely equivalent}; in
symbol: $X\approx_{\phi}Y$, if $H(X)=H(Y)$.
\end{definition}

We shall write%
\[
X^{\phi}=\{Y\in B:X\approx_{\phi}Y\}.
\]

It is clear that $X^{\phi}\subseteq X^{f}$. In the general case these classes
are not coincide. The set of all different classes of kind $X^{\phi}$ may be
regarded as the set of suitable subsets of $\frak{M}$ and will be denoted by
$\Phi(\mathcal{B})$.

\subsection{Topological properties of $\Phi(\mathcal{B})$}

Define on $\Phi(\mathcal{B})$ a topology $\mathcal{T}$, which base of open
sets consists of finite unions and intersections of sets of kind%
\[
V_{+}(A)=\{X^{\phi}:A\in H(X)\};\text{ \ }V_{-}(A)=\{X^{\phi}:A\notin H(X)\},
\]
where $A$ runs $\frak{M}$.

In the sequel there will be needed following notions.

Let $A_{1}$, $A_{2}$, ..., $A_{p}$ be finite dimensional Banach spaces. This
collection is said to be \textit{independent} if there is no pair $i$, $j$ for
which $A_{i}\in H(A_{j})$.

\begin{definition}
Let $X\in\mathcal{B}$. The set $\frak{M}(X^{f})$ has the \textit{compactness
property} (shortly: c. p.) if for any independent collection $A_{1}$, $A_{2}$,
.., $A_{n}$, $B_{1}$, $B_{2}$, ..., $B_{m}$ there exists $C\in\frak{M}%
(X^{f})$ such that $A_{1}$, $A_{2}$, ..., $A_{n}\in H(C)$ and $B_{1}$, $B_{2}%
$, ..., $B_{m}\notin H(C)$.
\end{definition}

Clearly, $\frak{M}(\left(  l_{\infty}\right)  ^{f})$ has the c. p. For any
finite dimensional $X\in\mathcal{B}$ the set $\frak{M}(X^{f})$ fails to have
this property.

\begin{theorem}
The space $\left\langle \Phi(\mathcal{B}),\mathcal{T}\right\rangle $ is
totally disconnected Hausdorff compact topological space.
\end{theorem}

\begin{proof}
$V_{+}(A)=\Phi(\mathcal{B})\backslash V_{-}(A)$ for every $A\in\frak{M}$.
Hence, every open set $U\in\mathcal{T}$ is also closed, i.e. $\left\langle
\Phi(B),\mathcal{T}\right\rangle $ is totally disconnected.

If $X^{\phi}\neq Y^{\phi}$ then $H(X)\neq H(Y)$, i.e. there exists
A$\in\frak{M}$ such that $A\in H(X)$ and $A\notin H(Y)$. It is clear that
$X^{\phi}\in V_{+}(A)$, $Y^{\phi}\in V_{-}(A)$ and, consequently,
$\left\langle \Phi(B),\mathcal{T}\right\rangle $ is Hausdorff.

To show that each centered system of closed sets of $\Phi(\mathcal{B})$ has a
nonempty intersection, using the c. p. and the Zorn lemma, extend a given
centered system $F$ of closed sets of $\Phi(\mathcal{B})$ to a maximal
centered system $F_{\max}$. It is clear that the intersection of $F_{\max}$
consists of the unique point - some class $X^{\phi}$, i.e. is not empty.
\end{proof}

\begin{remark}
The compactness property means that the sets $V_{+}(A)$ and $V_{-}(A)$
generate the Boolean algebra, say $\frak{B}^{\prime}$ with operations
\begin{align*}
&  V_{\pm}(A_{1})\wedge V_{\pm}(A_{2})\overset{\operatorname*{def}}{=}V_{\pm
}(A_{1})\cap V_{\pm}(A_{2}),\\
&  V_{\pm}(A_{1})\vee V_{\pm}(A_{2})\overset{\operatorname*{def}}{=}%
V_{+}(A_{1})\cup V_{\pm}(A_{2}),\\
&  \complement V_{+}(A_{1})\overset{\operatorname*{def}}{=}V_{-}(A_{1});\text{
}\complement V_{-}(A_{1})\overset{\operatorname*{def}}{=}V_{+}(A_{1}).
\end{align*}

Its Stonian space $\frak{S}\left(  \frak{B}^{\prime}\right)  $ is just the
topological space $\left\langle \Phi(\mathcal{B}),\mathcal{T}\right\rangle $,
as it follows from the definition.
\end{remark}

Any Hausdorff compact space is normal and the same is true for $\frak{M}$ when
the last set is regarded as the topological subspace of $\Phi(\mathcal{B})$
endowed with the induced topology (but not the metric $\varrho$!).

Hence the Wallman compact extension of $\frak{M}$ (which may be identified
with $\Phi(\mathcal{B})$ according to definitions) is homeomorphic to its
Stone-\v{C}ech compactification. From this observation the theorem follows

\begin{theorem}
$\operatorname{card}\Phi(\mathcal{B})=2^{\frak{c}}.$
\end{theorem}

\begin{proof}
From the c. p. and definitions it follows the existence of an infinite subset
of $\Phi(\mathcal{B})$ which is discrete in the induced topology. Thus, the
Stone-\v{C}ech compactification $\beta\Phi(\mathcal{B})$ contains the
Stone-\v{C}ech compactification of naturals $\beta\mathbb{N}$ and is of
cardinality $2^{\frak{c}}$ because of $\operatorname{card}\beta\mathbb{N}%
=2^{\frak{c}}$.
\end{proof}

\begin{remark}
It is unknown: whether every class $X^{\phi}$, which is generated by an
infinitely dimensional Banach space X, is a proper class (i.e. contains a
space of arbitrary large dimension).

Using some model theory it may be shown that there exists such cardinal number
$\tau$ that any Banach space of dimension $\geq\tau$ generates the proper
class $X^{\phi}$. As $\tau$ may be chosen so called the Hunf number
$h(\omega_{1},\omega)$ of the infinitary logic $\frak{L}_{\omega_{1},\omega}$.
It is known that $h(\omega_{1},\omega)=\gimel(\omega_{1})$.
\end{remark}

So, in the case $\dim(X)\geq\gimel(\omega_{1})$, the class $X^{\phi}$ contains
spaces of arbitrary dimension larger then $\dim\left(  X\right)  $.

\begin{remark}
Let $\frak{M}$ be regarded as the topological subspace of $\left\langle
\Phi(\mathcal{B}),\mathcal{T}\right\rangle $ (resp., of $\left\langle
f(\mathcal{B}),\mathcal{T}_{f}\right\rangle $ or $\left\langle \xi
(\mathcal{B}),\mathcal{T}_{\xi}\right\rangle $, equipped with the
corresponding restriction of topology ($\mathcal{T}\mid_{\frak{M}}$,
$\mathcal{T}_{f}\mid_{\frak{M}}$or $\mathcal{T}_{\xi}\mid_{\frak{M}}$
respectively). Then all these restrictions --$\mathcal{T}\mid_{\frak{M}}$,
$\mathcal{T}_{f}\mid_{\frak{M}}$or $\mathcal{T}_{\xi}\mid_{\frak{M}}$ are
coincide. The topological spaces $\left\langle \Phi(\mathcal{B}),\mathcal{T}%
\right\rangle $, $\left\langle f(\mathcal{B}),\mathcal{T}_{f}\right\rangle $
and $\left\langle \xi(\mathcal{B}),\mathcal{T}_{\xi}\right\rangle $ may be
regarded as corresponding compactifications of $\frak{M}$.

E.g., $\left\langle \Phi(\mathcal{B}),\mathcal{T}\right\rangle $ is just the
Stone-\v{C}ech compactification of $\left\langle \frak{M},\mathcal{T}%
\mid_{\frak{M}}\right\rangle .$

Notice that there exists a continuous surjection $h_{f}:\left\langle
\Phi(\mathcal{B}),\mathcal{T}\right\rangle \rightarrow\left\langle
f(\mathcal{B}),\mathcal{T}_{f}\right\rangle $ that holds $\frak{M}$ and glues
finitely equivalent elements of $\Phi(\mathcal{B})$.

Similarly, there exists a continuous surjection $h_{\xi}:\left\langle
\Phi(\mathcal{B}),\mathcal{T}\right\rangle \rightarrow\left\langle
\xi(\mathcal{B}),\mathcal{T}_{\xi}\right\rangle $ that holds $\frak{M}$ and
glues finitely equivalent elements of $\xi\mathcal{B})$.

Surely, from cardinality arguments it follows the existence of such different
classes $X^{\phi}$ and $Y^{\phi}$ that are finitely equivalent and, as well,
the existence of such different classes $X_{1}^{\phi}$ and $Y_{1}^{\phi}$ that
are elementary equivalent. It is easy to give examples of such classes.
\end{remark}

\begin{example}
Classes $(c_{0})^{\phi}$and $(l_{\infty})^{\phi}$ are finitely equivalent. At
the same time, for every finite-dimensional subspace of $c_{0}$ its unit ball
is a polyhedron. Certainly, $(c_{0})^{\phi}<_{\phi}$ $(l_{\infty})^{\phi}$ and
$(c_{0})^{\phi}$ is not strongly finitely equivalent to $(l_{\infty})^{\phi}$.
\end{example}

\begin{example}
Another example may be obtained if instead $(l_{\infty})^{\phi}$ it will be
regarded a class $(\widetilde{c_{0}})^{\phi}$, where $\widetilde{c_{0}}$ is
the space $c_{0}$, regarded under the equivalent strictly convex norm. Surely,
$(c_{0})^{\phi}$ and $(\widetilde{c_{0}})^{\phi}$ are finitely equivalent and
are strongly finitely incomparable.
\end{example}

\begin{example}
Classes $(l_{p})^{\phi}$ and $(l_{p}\oplus_{p}L_{p})^{\phi}$ are distinct
(since $l_{p}$ does not contain $l_{2}^{\left(  n\right)  }$ for large $n$).
Nevertheless, these classes are elementary equivalent (see e.g. \cite{hh})
\end{example}

\subsection{Omittable spaces in the class $X^{f}$}

Put%
\[
\Phi(X^{f})=\{Y^{\phi}:X\sim_{f}Y\}.
\]

\begin{definition}
A space $A\in\frak{M}(X^{f})$ is said to be \textit{omittable in the class
}$X^{f}$ if there exists such $Y^{\phi}\in\Phi(X^{f})$ that $A\notin H(Y)$ A
space $A\in\frak{M}(X^{f})$ is said to be \textit{non omittable} in the
contrary case.
\end{definition}

It is obvious that if $X$ is of finite dimension then every space $A\in H(X)$
is non omittable. Another example gives the class $\left(  l_{\infty}\right)
^{f}$, in which every nontrivial (i.e. more then one-dimensional) space from
$\frak{M}(\left(  l_{\infty}\right)  ^{f})$ is omittable.

\begin{theorem}
A space $A\in\frak{M}(X^{f})$ is non-omittable if and only if $A$ is an
isolated point in $\left\langle \Phi(\mathcal{B}),\mathcal{T}\right\rangle $.
\end{theorem}

\begin{proof}
If $X$ is finite dimensional then every $A\in\frak{M}(X^{f})$ is an isolated
point. Indeed, $V_{+}(X)=X$ and $X$ is isolated. If $A\in\frak{M}(X^{f})$ and
$\dim A=\dim X-1$ then $\{A\}=V_{+}(A)\cap V_{-}(X)$ and, hence, $A$ is an
isolated point. Then we are proceeding by induction. Let $X$ be infinite
dimensional. Let $\frak{A}$ be a set of all non omittable spaces in
$\frak{M}(X^{f})$. It is clear that $A$ satisfies properties (\textbf{H}) and
(\textbf{A}$_{0}$). From the theorem 3 follows that there exists a space
$Y\in\frak{B}$ with $H(Y)=\frak{A}$. If $Y$ is infinite dimensional then there
exists such $B\in H(Y)$ that $A\in H(B)$ and $B$ is non omittable. Then
$\{A\}=V_{+}(A)\cap V_{-}(B)$ is an isolated point. If $Y$ is of finite
dimension and $A\in\frak{A}$ then either $A{{}\neq} Y$ and is isolated as
below or $A=Y$. Consider a set of all $B\in\frak{M}(X^{f})$ such that $A\notin
H(B)$. This set does not have the property (\textbf{A}$_{0}$). Indeed, in the
contrary case the set $\{B\in\frak{M}(X^{f}):A\notin H\left(  B\right)  \}$
may be regarded as $H(Z)$ for some $Z\in\mathcal{B}$. However, this
contradicts with the non-omitability of $A$. Hence, there exists a space
$B\in\frak{M}(X^{f})$ such that $B$ contains a subspace, isometric to $A$. In
this case $A=V_{+}(A)\cap V_{-}(B)$ and, hence, is an isolated point too.
\end{proof}

\end{document}